\documentclass[twoside,leqno,10pt]{amsart}
\usepackage{amsfonts}
\usepackage{amsmath}
\usepackage{enumerate}
\usepackage{color}
\begin{document}

\newtheorem{theorem}{Theorem}[section]
\newtheorem{thm}{Theorem}
\newtheorem{lemma}{Lemma}[section]
\newtheorem{prop}{Proposition}[section]
\newtheorem{preim}{Lemma}[section]
\newtheorem{cor}{Corollary}[section]
\newtheorem{con}{Conjecture}[section]
\renewcommand{\thethm}{\Roman{thm}}
\numberwithin{equation}{section}
\newcommand{\QED}{{\hfill\large $\Box$}}

\theoremstyle{definition}
\newtheorem{definition}[theorem]{Definition}
\newtheorem{example}[theorem]{Example}
\newtheorem{xca}[theorem]{Exercise}

\theoremstyle{remark}
\newtheorem{remark}[theorem]{Remark}

\newcommand{\ps}{\pi/\sin\frac{\pi}{p}}
\newcommand{\nsum}{\sum^\infty_{n=1}}
\newcommand{\msum}{\sum^\infty_{m=1}}
\newcommand{\mmsum}{\sum^{k-1}_{m=1}}
\newcommand{\fint}{\int^\infty_0f(t)\,dt}
\newcommand{\fkint}{\int^\infty_kf(t)\,dt}
\newcommand{\rlam}{\lambda(r)}
\newcommand{\ffm}{\frac{1}{1+m}\left(\frac{1}{m}\right)^{\frac{1}{r}}}
\newcommand{\sumi}{\mathop{{\sum}^{'}}}

\title{On Copson's inequalities for $0<p<1$}
\author{Peng Gao and Huayu Zhao}
\subjclass[2000]{Primary 26D15} \keywords{Copson's inequalities}


\begin{abstract}
  Let $(\lambda_n)_{n \geq 1}$ be a non-negative sequence with
  $\lambda_1>0$ and let $\Lambda_n=\sum^n_{i=1}\lambda_i$. We study the following Copson inequality for $0<p<1$, $L>p$,
\begin{align*}
   \sum^{\infty}_{n=1}\left (\frac 1{\Lambda_n} \sum^{\infty}_{k=n}\lambda_k x_k \right
   )^p \geq \left ( \frac {p}{L-p}\right )^p
   \sum^{\infty}_{n=1}x^p_n.
\end{align*}
    We find conditions on $\lambda_n$ such that the above inequality is valid with the constant being best possible.
\end{abstract}

\maketitle
\section{Introduction}
\label{sec 1} \setcounter{equation}{0}

  Let $p>0$ and ${\bf x}=(x_n)_{n \geq 1}$ be a non-negative sequence.
  Let $(\lambda_n)_{n \geq 1}$ be a non-negative sequence with
  $\lambda_1>0$ and let $\Lambda_n=\sum^n_{i=1}\lambda_i$.
  The well-known Copson inequalities \cite[Theorem 1.1, 2.1]{C} state that:
\begin{align}
\label{1.1}
  & \sum^{\infty}_{n=1}\lambda_n\Lambda^{-c}_n\left ( \sum^{n}_{k=1}\lambda_k x_k \right
   )^p \leq \left ( \frac {p}{c-1}\right )^p
   \sum^{\infty}_{n=1}\lambda_n\Lambda^{p-c}_nx^p_n, \ 1<c \leq p
   ; \\
  \label{1.2}
  & \sum^{\infty}_{n=1}\lambda_n\Lambda^{-c}_n\left ( \sum^{\infty}_{k=n}\lambda_k x_k \right
   )^p \leq \left ( \frac {p}{1-c}\right )^p
   \sum^{\infty}_{n=1}\lambda_n\Lambda^{p-c}_nx^p_n, \ 0 \leq c < 1<p.
\end{align}

   The above two inequalities are equivalent (see \cite{G10}) and the constants are best possible. When $\lambda_k=1, k \geq 1$ and $c=p$, inequality \eqref{1.1} becomes the following celebrated
   Hardy inequality (\cite[Theorem 326]{HLP}):
\begin{align}
\label{1.3}
  \sum^{\infty}_{n=1}\Big( \frac 1{n} \sum^{n}_{k=1}x_k \Big
  )^p \leq \left ( \frac {p}{p-1}\right )^p \sum^{\infty}_{n=1}x^p_n.
\end{align}

    Note that the reversed inequality of \eqref{1.2} holds when $c \leq 0 <p <1$ (\cite[Theorem 2.3]{C}) with the constant being best possible and as pointed out in \cite[p. 390]{B2}, the reversed inequality of \eqref{1.2} continues to hold with constant $p^p$ when $c>0$. The particular case of $c=p, \lambda_k=1, k \geq 1$ becomes the following one given in \cite[ Theorem 345]{HLP}:
\begin{align}
\label{1.4}
  \sum^{\infty}_{n=1}\Big( \frac 1{n} \sum^{\infty}_{k=n}x_k \Big
  )^p \geq p^p \sum^{\infty}_{n=1}x^p_n.
\end{align}

    It is noted in \cite{HLP} that the constant $p^p$ in \eqref{1.4} may not be best possible
  and the best constant for $0<p \leq 1/3$ was shown by Levin and Ste\v ckin
  \cite[Theorem 61]{L&S} to be indeed  $(p/(1-p))^p$. In \cite{G9}, it is shown that the constant $(p/(1-p))^p$ stays best possible for all $0<p \leq 0.346$.
  It is further shown in \cite{G11} that the constant $(p/(1-p))^p$ is best possible when $p=0.35$.

  There exists an extensive and rich literature on extensions and generalizations of Copson's inequalities and Hardy's inequality \eqref{1.3} for $p>1$. For recent developments in this direction, we refer the reader to the articles in \cite{G6}-\cite{G11} and the references therein. On the contrary, the case $0<p<1$ is less known as this can be seen by comparing inequalities \eqref{1.3} and \eqref{1.4}. On one hand, the constant in \eqref{1.3} is shown to be best possible for all $p>1$. On the other hand, though it is known the best constant that makes inequality \eqref{1.4} valid is $(p/(1-p))^p$ when $0<p \leq 0.346$, it is shown in \cite{G9} that the constant $(p/(1-p))^p$ fails to be best possible when $1/2 \leq p <1$ and the best constant in these cases remains unknown.

  Our goal in this paper is to study the following variation of Copson's inequalities for $0<p<1$:
\begin{align}
\label{1.5}
   \sum^{\infty}_{n=1}\left (\frac 1{\Lambda_n} \sum^{\infty}_{k=n}\lambda_k x_k \right
   )^p \geq \left ( \frac {p}{L-p}\right )^p
   \sum^{\infty}_{n=1}x^p_n,
\end{align}
    where $L>p$ is a constant.

    Let $q<0$ be the number satisfying $1/p+1/q=1$, then inequality \eqref{1.5} is equivalent to its dual version (assuming that $x_n>0$ for all $n$):
\begin{align}
\label{1.5'}
   \sum^{\infty}_{n=1}\left (\lambda_n  \sum^{n}_{k=1} \frac {x_k}{\Lambda_k} \right
   )^q \leq \left ( \frac {p}{L-p}\right )^q
   \sum^{\infty}_{n=1}x^q_n.
\end{align}
   The equivalence of the above two inequalities can be easily established following the discussions in \cite[Section 1]{G9}.

    Our main result gives a condition on $\lambda_n$ and $L$ such that inequalities \eqref{1.5} and \eqref{1.5'} hold. For this purpose, we define
for constants $p$ and $L$,
\begin{align}
\label{1.8}
   a_1(L,p) =& (\frac {L}p-2)^2(1+L\frac {2-p}{1-p}) \\
   &-(1+(\frac Lp-2)\frac {1-2p}{1-p})(L^2(L-1)^2+2L(L-1)(L-p-1)+L^2-2(L-1)(p+1)),  \nonumber \\
   a_2(L,p) =& (\frac 1{p}-1)L^4+\frac {(1-p)(1-2p)}{p}L^3-(3-p)(1-p)L^2-(p^2-p+2)L+2p(1+p). \nonumber
\end{align}

   Our main result is the following
\begin{theorem}
\label{thm1}
    Let $0<p<1$ be fixed.  Let $(\lambda_n)_{n \geq 1}$ be a non-negative sequence with
  $\lambda_1>0$ and let $\Lambda_n=\sum^n_{i=1}\lambda_i$. If there exists a positive constant
    $L>p$ such that for any integer $n \geq 1$,
\begin{align}
\label{1.6}
  \frac {L-p}p \cdot \frac {\lambda_n}{\Lambda_n} \leq (1+(\frac Lp-2)\frac {\lambda_n}{\Lambda_n})^{1/(1-p)}-(1-\frac {\lambda_{n+1}}{\Lambda_{n+1}})^{(1+p)/(1-p)}(\frac {\lambda_n}{\Lambda_{n}})^{1/(1-p)}(\frac {\lambda_{n+1}}{\Lambda_{n+1}})^{-1/(1-p)},
\end{align}
    then inequality \eqref{1.5} holds for all non-negative sequences ${\bf x}$. In particular, suppose that
\begin{align}
\label{1.7}
    L=\sup_n\Big(\frac {\Lambda_{n+1}}{\lambda_{n+1}}-\frac
    {\Lambda_n}{\lambda_n}\Big) > p ~~.
\end{align}
   Let $a_i(L,p), i=1,2$ be defined as in \eqref{1.8}, then inequality \eqref{1.5} holds for all non-negative sequences ${\bf x}$ when $L \geq 1,0<p \leq 1/3$
   and $a_1(L,p)\geq 0$ or when $0<L < 1,0<p \leq L/4$ and $a_2(L,p) \geq 0$.
\end{theorem}

    We note that as $\lim_{L \rightarrow \infty}a_1(L,p)<0, \lim_{p \rightarrow 0^+}a_2(4p,p)/p<0$, the values of $a_i(L,p), i=1,2$ do give
    restrictions on the validity of inequality \eqref{1.5}.
    We note that when $\lambda_n=1$, then $L=1$ and Theorem \ref{thm1} implies the above mentioned result of Levin and Ste\v ckin. Here the choice of our condition on $L$ in \eqref{1.7} is based on the study of the $p>1$ analogue of inequality \eqref{1.5}, in which case a result of Cartlidge \cite{Car} shows that when the reversed inequality \eqref{1.7} holds for $p>1$, then the following inequality holds for all non-negative sequences ${\bf x}$:
\begin{align}
\label{1.9}
   \sum^{\infty}_{n=1}\left (\frac 1{\Lambda_n} \sum^{n}_{k=1}\lambda_k x_k \right
   )^p \leq \left ( \frac {p}{p-L}\right )^p
   \sum^{\infty}_{n=1}x^p_n, \quad p>1.
\end{align}

   In \cite{Be1}, \cite{G}-\cite{G9}, two special cases of inequality \eqref{1.9} corresponding to $\lambda_n=n^{\alpha}-(n-1)^{\alpha}$ and $\lambda_n=n^{\alpha-1}$ for $p>1, \alpha p >1$ were studied. It follows from these work that we know inequality \eqref{1.9} holds in either case with best possible constant $(\alpha p/(\alpha p-1))^p$ except for the case when $\lambda_n=n^{\alpha-1}$, $1<p \leq 4/3, 1\leq \alpha <1+1/p$ or $4/3<p<2, 1 \leq \alpha<2$.

  It is now interesting to study the following $0<p<1$ analogues of the above inequalities:
\begin{align}
\label{1.12}
   \sum^{\infty}_{n=1}\Big (\frac
1{n^{\alpha}}\sum^{\infty}_{k=n}(k^{\alpha}-(k-1)^{\alpha})x_k\Big )^p &
\geq  \Big( \frac {\alpha p}{1-\alpha p} \Big )^p\sum^{\infty}_{n=1}x_n^p, \\
\label{1.13}
   \sum^{\infty}_{n=1}\Big (\frac
1{\sum^n_{i=1}i^{\alpha-1}}\sum^{\infty}_{k=n}k^{\alpha-1}x_k\Big )^p &
\leq  \Big(\frac {\alpha p}{1-\alpha p} \Big
)^p\sum^{\infty}_{n=1}x_n^p.
\end{align}

    Note that when $\alpha=1$, the above two inequalities become inequality \eqref{1.4}. We note that it is shown in \cite[(4.14)]{G9}
    that inequality \eqref{1.13} holds when $0< \alpha <1, 0<p<1/(\alpha+2)$. It is also shown in \cite[Theorem 1.1]{G8} that inequality \eqref{1.12} holds when $\alpha>0, 0<p<1/(\alpha+2)$ when one replaces $k^{\alpha}-(k-1)^{\alpha}$ by $(k+1)^{\alpha}-k^{\alpha}$. As $(k+1)^{\alpha}-k^{\alpha} \leq k^{\alpha}-(k-1)^{\alpha}$ when $0<\alpha \leq 1$, this implies that inequality \eqref{1.12} holds when $0<\alpha \leq 1, 0<p<1/(\alpha+2)$.

   Note that the values of $p$ are not given explicitly in Theorem \ref{thm1} when $0<L<1$. Moreover, the condition \eqref{1.7}
   may not always be satisfied. For these reasons, and with future applications in mind, we develop the following
\begin{theorem}
\label{thm1'}
    Let $0<p<1$ be fixed.  Let $(\lambda_n)_{n \geq 1}$ be a non-negative sequence with
  $\lambda_1>0$ and let $\Lambda_n=\sum^n_{i=1}\lambda_i$. If \eqref{1.7} is satisfied with $0<L<1$ and that
\begin{align}
\label{1.10}
  p \leq \frac {L^2}{4}:=p_L,
\end{align}
   then inequality \eqref{1.5} holds for all non-negative sequences ${\bf x}$.

   If there exist positive constants
    $1/2< L<1, 0<M<1, L+2M<1$ such that for any integer $n \geq 1$,
\begin{align}
\label{1.15}
    \frac{\Lambda_{n+1}}{\lambda_{n+1}}-\frac
    {\Lambda_n}{\lambda_n} \leq L+M\frac {\lambda_n}{\Lambda_n}.
\end{align}
    Then inequality \eqref{1.5} holds for all non-negative sequences ${\bf x}$ when
\begin{align}
\label{1.16}
  p \leq \min \Big \{ \frac {L(2L-1)}{4(2L+M)}, \frac {L(1-L-2M)}{2(1-L-M)} \Big \}.
\end{align}
\end{theorem}

   We remark here that it is easy to see that the minimum on the right-hand side of \eqref{1.16} can take either values.  We now apply Theorem \ref{thm1'} to study inequalities \eqref{1.12}-\eqref{1.13}.  As the case $\alpha=1$ yields the classical inequality \eqref{1.4}, we concentrate on the case $\alpha>1$ and we deduce readily from Theorem \ref{thm1'} the following
\begin{theorem}
\label{thm2}
    Let $\alpha \geq 1$ and $p_{1/\alpha}$ be defined as in \eqref{1.10}. Then inequality \eqref{1.12} holds for all non-negative sequences ${\bf x}$ when $\alpha > 1, 0< p \leq p_{1/\alpha}$ and inequality \eqref{1.13} holds for all non-negative sequences ${\bf x}$ when $\alpha \geq 2, 0< p \leq p_{1/\alpha}$. The constants are best possible in both cases.
\end{theorem}

   In fact, note that \cite[Lemma 2.1]{G} implies that \eqref{1.7} is satisfied for $\lambda_n=n^{\alpha}-(n-1)^{\alpha}$ with $L=1/\alpha$ when $\alpha \geq 1$ and \eqref{1.7} is satisfied for $\lambda_n=n^{\alpha-1}$ with $L=1/\alpha$ when $\alpha \geq 2$. That the constant is best possible can be seen by setting
   $x_n=n^{-1/p-\epsilon}$ and letting $\epsilon \rightarrow 0^+$.

\section{Proof of Theorem \ref{thm1}}
\label{sec 2} \setcounter{equation}{0}
   Our goal is to find conditions on the $\lambda_n$ such that the following inequality holds for $0<p<1$, $L>p$:
\begin{align*}
   \sum^{\infty}_{n=1}\left ( \frac 1{\Lambda_n}\sum^{\infty}_{k=n}\lambda_k x_k \right
   )^p \geq \left ( \frac {p}{L-p}\right )^p
   \sum^{\infty}_{n=1}x^p_n.
\end{align*}
   It suffices to prove the above inequality by replacing the infinite sums by finite sums from $n=1$ to $N$ (and $k=n$ to $N$) for any integer $N \geq 1$. Note that as in \cite[Section 3]{G6}, we have
\begin{align*}
   \sum^N_{n=1}a^p_n =\sum^N_{n=1} \frac {a^p_n}{\sum^n_{i=1} w_i}\sum^n_{k=1} w_k=\sum^N_{n=1} w_n\sum^N_{k=n} \frac {a^p_k}{\sum^k_{i=1} w_i}.
\end{align*}
   By H\"older's inequality, we have
\begin{align*}
  \sum^N_{k=n} \frac {a^p_k}{\sum^k_{i=1} w_i} \leq \Big ( \sum^N_{k=n} \Big (\lambda^p_k \sum^k_{i=1} w_i\Big )^{-1/(1-p)}\Big )^{1-p}\Big ( \sum^N_{k=n}\lambda_k a_k \Big )^{p}.
\end{align*}
   It follows that
\begin{align*}
   \sum^N_{n=1}a^p_n  \leq \sum^N_{n=1} w_n\Big ( \sum^N_{k=n} \Big (\lambda^p_k \sum^k_{i=1} w_i\Big )^{-1/(1-p)}\Big )^{1-p}\Big ( \sum^N_{k=n}\lambda_k a_k \Big )^{p}.
\end{align*}

   Suppose we can find a sequence ${\bf w}$ of positive terms, such that for any integer $n\geq 1$,
\begin{align*}
  (\sum^n_{i=1}w_i)^{1/(p-1)} \leq \Big (\frac p{L-p} \Big )^{p/(p-1)}\lambda^{p/(1-p)}_n\Big (\frac {w^{1/(p-1)}_n}{\Lambda^{p/(1-p)}_n}-\frac {w^{1/(p-1)}_{n+1}}{\Lambda^{p/(1-p)}_{n+1}} \Big ).
\end{align*}
   Then the desired inequality follows. Now we make a change of variables: $w_i \mapsto \lambda_iw_i$ to recast the above inequality as
\begin{align}
\label{2.1}
  (\sum^n_{i=1}\lambda_iw_i)^{1/(p-1)} \leq \Big (\frac p{L-p} \Big )^{p/(p-1)}\lambda^{p/(1-p)}_n\Big (\frac {\lambda^{1/(p-1)}_nw^{1/(p-1)}_n}{\Lambda^{p/(1-p)}_n}-\frac {\lambda^{1/(p-1)}_{n+1}w^{1/(p-1)}_{n+1}}{\Lambda^{p/(1-p)}_{n+1}} \Big ).
\end{align}
   We now define our sequence ${\bf w}$ to satisfy: $w_1=1$ and inductively that for $n \geq 2$,
\begin{align}
\label{2.2}
  \frac 1{\Lambda_n}\sum^n_{i=1}\lambda_iw_i=\frac {p}{L-p}w_{n+1}.
\end{align}
   This allows us to deduce that for $n \geq 1$,
\begin{align}
\label{2.3}
  w_{n+1}=(1+(\frac Lp-2)\frac {\lambda_n}{\Lambda_n})w_n.
\end{align}
   Applying \eqref{2.2}, \eqref{2.3} in \eqref{2.1}, we see that inequality \eqref{2.1} becomes \eqref{1.6}.

   We now set $x=\lambda_n/\Lambda_n, y=\lambda_{n+1}/\Lambda_{n+1}$ to recast inequality \eqref{1.6} as
\begin{align*}
  \frac {L-p}p x \leq (1+(\frac Lp-2)x)^{1/(1-p)}-(\frac 1y-1)^{(1+p)/(1-p)}x^{1/(1-p)}y^{p/(1-p)}.
\end{align*}

   As the function $t \mapsto (t-1)^{(1+p)/(1-p)}t^{-p/(1-p)}$ is an increasing function of $t \geq 1$, and that condition \eqref{1.7} implies that
   $1/y \leq 1/x+L$, it suffices to establish the above inequality for $1/y=1/x+L$. To facilitate the proof of Theorem \ref{thm1'}, we proceed by
  taking the weaker condition \eqref{1.15} into consideration to assume
   that $1/y \leq 1/x+L+Mx$, where $M \geq 0$ is a constant.  Again in this case,
   it suffices to prove the above inequality for $1/y=1/x+L+Mx$, which is equivalent to showing that $f_{L,M,p}(x) \geq 0$, where
\begin{align*}
   & f_{L,M,p}(x) \\
   := & (1+(\frac Lp-2)x)^{1/(1-p)}-(1+(L-1)x+Mx^2)^{(1+p)/(1-p)}(1+Lx+Mx^2)^{-p/(1-p)}-\frac {L-p}p x.
\end{align*}
  Suppose that $L \geq 1$ and \eqref{1.7} is valid. In this case we can set $M=0$ so that it suffices to show that $f_{L,0,p}(x) \geq 0$. Calculation shows that
\begin{align*}
  \frac {(1-p)^2}{p}f''_{L,0,p}(x)=(1+(\frac Lp-2)x)^{(2p-1)/(1-p)}(1+(L-1)x)^{2p/(1-p)-1}(1+Lx)^{-1/(1-p)-1}g_{L,p}(x),
\end{align*}
   where
\begin{align*}
  & g_{L,p}(x)  \\
  = &(\frac {L}p-2)^2(1+(L-1)x)^{\frac {1-3p}{1-p}}(1+Lx)^{\frac {2-p}{1-p}}\\
  &-(1+(\frac Lp-2)x)^{\frac {1-2p}{1-p}}(L^2(L-1)^2x^2+2L(L-1)(L-p-1)x+L^2-2(L-1)(p+1)).
\end{align*}

   Suppose that $0<p \leq 1/3$. We want to show that $g_{L,p}(x) \geq 0$ for $0<x \leq 1$. We first note that we have
\begin{align*}
   g_{L,p}(x) \geq &(\frac {L}p-2)^2(1+(L-1)x)^{\frac {1-3p}{1-p}}(1+Lx)^{\frac {2-p}{1-p}}\\
  &-(1+(\frac Lp-2)x)^{\frac {1-2p}{1-p}}(L^2(L-1)^2x+2L(L-1)(L-p-1)x+L^2-2(L-1)(p+1)).
\end{align*}

   We may now assume that
\begin{align*}
    L^2(L-1)^2x+2L(L-1)(L-p-1)x+L^2-2(L-1)(p+1) \geq 0.
\end{align*}
   For otherwise, we have trivially $g_{L,p}(x) \geq 0$.
   We then estimate $(1+(L-1)x)^{\frac {1-3p}{1-p}}$ trivially by $(1+(L-1)x)^{\frac {1-3p}{1-p}} \geq 1$ and we apply Taylor expansion to see that
\begin{align*}
  \quad (1+Lx)^{\frac {2-p}{1-p}} & \geq 1+L\frac {2-p}{1-p}x, \\
  (1+(\frac Lp-2)x)^{\frac {1-2p}{1-p}} & \leq 1+(\frac Lp-2)\frac {1-2p}{1-p}x.
\end{align*}

   It then follows that $g_{L,p}(x) \geq u_{L,p}(x)$, where
\begin{align*}
  u_{L,p}(x) = & (\frac {L}p-2)^2(1+L\frac {2-p}{1-p}x) \\
  &-(1+(\frac Lp-2)\frac {1-2p}{1-p}x)(L^2(L-1)^2x+2L(L-1)(L-p-1)x+L^2-2(L-1)(p+1)).
\end{align*}
   It is easy to see that $u_{L,p}(x)$ is a concave function, hence is minimized at $x=0$ or $x=1$. One checks that in this case
\begin{align}
\label{2.5}
  u_{L,p}(0)=(\frac {L}p-2)^2  - L^2+2(L-1)(p+1).
\end{align}
   When we regard the above expression as a function of $L$, it is easy to see that $u_{L,p}(0) \geq 0$ when $p \leq 1/3$. On the other hand, we have
   $u_{L,p}(1)=a_1(L,p) \geq 0$ by our assumption,  where $a_1(L,p)$ is defined in \eqref{1.8}. It follows that
   $g_{L,p}(x) \geq u_{L,p}(x) \geq 0$, hence $f''_{L,0,p}(x) \geq 0$. As $f_{L,0,p}(0)=f'_{L,0,p}(0)=0$, we then deduce that
   $f_{L,0,p}(x) \geq 0$ and this completes the proof for the case $L \geq 1$.

    When $0<L<1$, we note that
\begin{align*}
   &\frac p{L-p}f'_{L,M,p}(x)+1 \\
   =& (\frac Lp-2)\frac p{(1-p)(L-p)}(1+(\frac Lp-2)x)^{p/(1-p)} \\
  &  +\frac {p(1+p)(1-L)}{(1-p)(L-p)}(1-\frac {2Mx}{1-L})(1+(L-1)x+Mx^2)^{(2p)/(1-p)}(1+Lx+Mx^2)^{-p/(1-p)} \\
  &+\frac {Lp^2}{(1-p)(L-p)}(1+\frac {2Mx}{L})(1+(L-1)x+Mx^2)^{(1+p)/(1-p)}(1+Lx+Mx^2)^{-1/(1-p)}.
\end{align*}

   One checks that
\begin{align*}
  (\frac Lp-2)\frac p{(1-p)(L-p)}+\frac {p(1+p)(1-L)}{(1-p)(L-p)}+\frac {Lp^2}{(1-p)(L-p)}=1.
\end{align*}
   It follows from this and the arithmetic-geometric mean inequality that
\begin{align}
\label{f'}
  & \frac p{L-p}f'_{L,M,p}(x)+1 \\
  \geq &(1+(\frac Lp-2)x)^{\frac p{1-p}\cdot(\frac Lp-2)\frac p{(1-p)(L-p)}} \nonumber \\
  &  \cdot \Big ( (1-\frac {2Mx}{1-L})(1+(L-1)x+Mx^2)^{(2p)/(1-p)}(1+Lx+Mx^2)^{-p/(1-p)} \Big )^{\frac {p(1+p)(1-L)}{(1-p)(L-p)}} \nonumber \\
  & \cdot \Big ((1+\frac {2Mx}{L})(1+(L-1)x+Mx^2)^{(1+p)/(1-p)}(1+Lx+Mx^2)^{-1/(1-p)} \Big )^{\frac {Lp^2}{(1-p)(L-p)}}. \nonumber
\end{align}

   Thus, it suffices to show the right-hand side expression above is $\geq 1$. We now assume that \eqref{1.7} is valid so that we can set $M=0$ in the above expression
   to see that it is equivalent to showing that
\begin{align*}
  h_{L,p}(x):=(1+(\frac Lp-2)x)^{L-2p}(1+(L-1)x)^{p(1+p)(2-L)}(1+Lx)^{-p(1+p-pL)} \geq 1.
\end{align*}
   Calculation shows that
\begin{align*}
  h'_{L,p}(x) =(1+(\frac Lp-2)x)^{L-2p-1}(1+(L-1)x)^{p(1+p)(2-L)-1}(1+Lx)^{-p(1+p-pL)-1}v_{L,p}(x),
\end{align*}
   where
\begin{align*}
  v_{L,p}(x) =&(L-2p)(\frac Lp-2)(1+(L-1)x)(1+Lx)-p(1+p)(2-L)(1-L)(1+(\frac Lp-2)x)(1+Lx) \\
  &-p(1+p-pL)L(1+(\frac Lp-2)x)(1+(L-1)x).
\end{align*}
   It is easy to see that $v_{L,p}(x)$ is a concave function of $x$ when $L \geq 2p$. It follows that $v_{L,p}(x) \geq \min \{ v_{L,p}(0), v_{L,p}(1) \}$ and
  one checks that $v_{L,p}(0)=p u_{L,p}(0)$, where $u_{L,p}(0)$ is defined in \eqref{2.5}. Similar to our discussions for the case $L>1$, one checks
  that $u_{L,p}(0) \geq 0$ when $L \geq 4p$. On the other hand, we have $v_{L,p}(1)=a_2(L,p) \geq 0$  by our assumption,
  where $a_2 (L,p)$ is defined in \eqref{1.8}. It follows that
   $h'_{L,p}(x) \geq 0$, hence $h_{L,p}(x) \geq h_{L,p}(0) \geq 1$ and this completes the proof for the case $0<L<1$.

\section{Proof of Theorem \ref{thm1'}}
\label{sec 3} \setcounter{equation}{0}

   First, we assume that \eqref{1.7} is valid and in this case, it suffices to find the value of $p$ such that $a_2(L,p) \geq 0$ by Theorem \ref{thm1}.
   Note that $\lim_{p \rightarrow 0^+}a_2(a p,p)/p<0$ when $a >1$, it is therefore not possible to show $a_2(L,p) \geq 0$ by assuming that $p \leq L/a$ for any
   $a>1$. We therefore seek to show $a_2(L,p) \geq 0$ for $p \leq L^2/4$.
   We first note that
\begin{align*}
  \frac {a_2(L,p)}{1-p} &\geq  \frac 1{p}L^4+\frac {1-2p}{p}L^3-(3-p)L^2+(p-\frac {2}{1-p})L+\frac {2p(1+p)}{1-p} \\
  & \geq 4L^2+4(1-2p)L-(3-p)L^2+(p-\frac {2}{1-p})L+\frac {2p(1+p)}{1-p}L \\
  & \geq 4(1-2p)L+(p-\frac {2}{1-p})L+2pL \\
  &=(4-5p-\frac {2}{1-p})L \geq 0,
\end{align*}
   where the last inequality follows from the observation that the function $p \mapsto 4-5p-\frac {2}{1-p}$ is non-negative for $0< p \leq 1/4$.
 This completes the proof for the first assertion of Theorem \ref{thm1'}.

   To prove the second assertion of Theorem \ref{thm1'}, we see from the proof of Theorem \ref{thm1} that it suffices to show the right-hand side expression of \eqref{f'} is $\geq 1$. We simplify it to see that it is equivalent to showing that
\begin{align}
\label{3.1}
  & (1+(\frac Lp-2)x)^{L-2p}(1-\frac {2Mx}{1-L})^{(1-p^2)(1-L)} \\
  &\cdot (1+\frac {2Mx}{L})^{p(1-p)L}(1+(L-1)x+Mx^2)^{p(1+p)(2-L)}(1+Lx+Mx^2)^{-p(1+p-pL)} \geq 1. \nonumber
\end{align}
   We assume that
\begin{align}
\label{3.02}
  \frac Lp-2 \geq 0.
\end{align}
    This implies that the function
\begin{align*}
 (1+(\frac Lp-2)x)(1-\frac {2Mx}{1-L})
\end{align*}
   is a concave function of $x$, hence is minimized at $x=0$ or $1$. When $x=0$, the above function takes value $1$.
We further assume that the above function takes value $\geq 1$ when $x=1$. That is,
\begin{align}
\label{3.2'}
 (\frac Lp-1)(1-\frac {2M}{1-L}) \geq 1.
\end{align}
   We then deduce that
\begin{align*}
 1-\frac {2Mx}{1-L} \geq (1+(\frac Lp-2)x)^{-1}.
\end{align*}
   We apply the above estimation and the estimation that $1+2Mx/L \geq 1$ in \eqref{3.1} to see that it suffices to show
\begin{align*}
  & h_{L,M,p}(x)\\
  &:= (1+(\frac Lp-2)x)^{L-2p-(1-p^2)(1-L)}(1+(L-1)x+Mx^2)^{p(1+p)(2-L)}(1+Lx+Mx^2)^{-p(1+p-pL)} \geq 1.
\end{align*}
   We now assume that
\begin{align}
\label{3.04}
   L-2p-(1-p^2)(1-L) >0.
\end{align}

  Then calculation shows that
\begin{align*}
  h'_{L,M,p}(x)&:=(1+(\frac Lp-2)x)^{L-2p-(1-p^2)(1-L)-1}(1+(L-1)x+Mx^2)^{p(1+p)(2-L)-1} \\
  & \cdot(1+Lx+Mx^2)^{-p(1+p-pL)-1}u_{L,M,p}(x),
\end{align*}
   where
\begin{align*}
  u_{L,M,p}(x) =&( L-2p-(1-p^2)(1-L))(\frac Lp-2)(1+(L-1)x+Mx^2)(1+Lx+Mx^2) \\
  &-p(1+p)(2-L)(1-L)(1+(\frac Lp-2)x)(1-\frac {2Mx}{1-L})(1+Lx +Mx^2) \\
  &-p(1+p-pL)L(1+(\frac Lp-2)x)(1+\frac {2Mx}{L})(1+(L-1)x +Mx^2) \\
  \geq &( L-2p-(1-p^2)(1-L))(\frac Lp-2)(1+(L-1)x )(1+Lx ) \\
  &-p(1+p)(2-L)(1-L)(1+(\frac Lp-2)x)(1+(L+M)x) \\
  &-p(1+p-pL)L(1+(\frac Lp-2)x)(1+\frac {2M}{L})(1+(L+M-1)x) :=v_{L,M,p}(x).
\end{align*}
   It is easy to see that $v_{L,M,p}(x)$ is a concave function of $x$ when
\begin{align*}
  &( L-2p-(1-p^2)(1-L))(\frac Lp-2)(L-1)L-p(1+p)(2-L)(1-L)(\frac Lp-2)(L+M) \\
  &-p(1+p-pL)(\frac Lp-2)(L+2M)(L+M-1)\leq 0.
\end{align*}
   We can recast the above inequality as
\begin{align}
\label{3.2}
  &( L-2p-(1-p^2)(1-L))(1-L)L+p(1+p)(2-L)(1-L)(L+M) \\
  &-p(1+p-pL)(L+2M)(1-L-M)\geq 0. \nonumber
\end{align}
  Assuming the above inequality, we see that $v_{L,M,p}(x) \geq \min \{ v_{L,M,p}(0), v_{L,M,p}(1) \}$ and we have
\begin{align*}
  v_{L,M,p}(0) =&( L-2p-(1-p^2)(1-L))(\frac Lp-2)-p(1+p)(2-L)(1-L)-p(1+p-pL)(L+2M), \\
   v_{L,M,p}(1) =&( L-2p-(1-p^2)(1-L))(\frac Lp-2)L(L+1)-p(1+p)(2-L)(1-L)(\frac Lp-1)(1+L+M) \\
   &-p(1+p-pL)(\frac Lp-1)(L+2M)(L+M).
\end{align*}

   We note that when $M < 1$,
\begin{align*}
    &p(1+p)(2-L)(1-L)+p(1+p-pL)(L+2M) \leq 2p(1+p)(1+L+M),\\
  & p(1+p)(2-L)(1-L)(\frac Lp-1)(1+L+M)+p(1+p-pL)(\frac Lp-1)(L+2M)(L+M) \\
  \leq & (1+p)L(1+L+M)+(1+p)L(L+2M)(L+M) \\
  \leq & (1+p)L(1+L+M)+(1+p)L(1+L+M)(L+M)= (1+p)L(1+L+M)^2.
\end{align*}

    We then deduce that $v_{L,M,p}(0) \geq 0$ when
\begin{align}
\label{3.3}
    ( L-2p-(1-p^2)(1-L))(\frac Lp-2) \geq 2p(1+p)(1+L+M),
\end{align}
    and that $v_{L,M,p}(1) \geq 0$ when
\begin{align}
\label{3.4}
    ( L-2p-(1-p^2)(1-L))(\frac Lp-2) \geq 2(1+p)(1+L+M).
\end{align}

   Thus, one just needs to find values of $p$ to satisfy inequalities \eqref{3.02}-\eqref{3.4}. We first note that
\begin{align*}
  & p(1+p)(2-L)(1-L)(L+M) \geq 0, \\
  & p(1+p-pL)(L+2M)(1-L-M) \leq p(1+p)(1-L)(1+L+M).
\end{align*}
   We apply the above estimates to see that inequality \eqref{3.2} is a consequence of the following inequality:
\begin{align}
\label{3.5}
    ( L-2p-(1-p^2)(1-L))\frac {L}{p} \geq (1+p)(1+L+M).
\end{align}

   As inequality \eqref{3.4} implies inequalities \eqref{3.3} and \eqref{3.5}, we first find values of $p$ so that inequality \eqref{3.4} holds. To do so, we first simply inequality \eqref{3.4} by noting that
\begin{align*}
   L-2p-(1-p^2)(1-L) \geq 2L-1-2p.
\end{align*}

    Using this, we see that it suffices to find values of $p$ to satisfy
\begin{align*}
    ( 2L-1-2p)(\frac Lp-2) \geq 2(1+p)(1+L+M).
\end{align*}

    We recast the above  as
\begin{align*}
   (2L-1)L-2(4L+M)p+2(1-L-M)p^2 \geq 0.
\end{align*}
   Note as $1-L-M>0$, we have $(1-L-M)p^2 \geq 0$, so that the above inequality follows from
\begin{align*}
   (2L-1)L-4(2L+M)p \geq 0,
\end{align*}
   which implies that
\begin{align}
\label{3.9}
   p \leq \frac {L(2L-1)}{4(2L+M)}.
\end{align}

    One checks that the above inequality implies inequalities \eqref{3.02} and \eqref{3.04}.  One further notes that inequality \eqref{3.2'} is equivalent to
\begin{align}
\label{3.10}
   p \leq \frac {L(1-L-2M)}{2(1-L-M)}.
\end{align}

   Combining inequalities \eqref{3.9} and \eqref{3.10}, one readily deduces the second assertion of Theorem \ref{thm1'} and this completes the proof of Theorem \ref{thm1'}.
\newline
\noindent{\bf Acknowledgments.} P. G. is supported in part by NSFC grant 11371043.


\vspace*{.5cm}
\noindent\begin{tabular}{p{8cm}p{8cm}}
School of Mathematics and Systems Science & Academy of Mathematics and Systems Science \\
Beihang University &  Chinese Academy of Sciences\\
Beijing 100191, China & Beijing 100190, China \\
Email: {\tt penggao@buaa.edu.cn} & Email: {\tt zhaohuayu17@mails.ucas.ac.cn} \\
\end{tabular}
\end{document}